\documentclass[preprint,showpacs,preprintnumbers,amsmath,amssymb]{revtex4}

\usepackage{graphicx}
\usepackage{dcolumn}
\usepackage{bm}

\newcommand{\ba}[1]{\begin{array}{#1}}
\newcommand{\ea}{\end{array}}

\newcommand{\be}{\begin{equation}}
\newcommand{\ee}{\end{equation}}
\newcommand{\bea}{\begin{eqnarray}}
\newcommand{\eea}{\end{eqnarray}}
\newcommand{\beann}{\begin{eqnarray*}}
\newcommand{\eeann}{\end{eqnarray*}}

\begin{document}

\title{Monte Carlo Tests of SLE Predictions for the 2D Self-Avoiding Walk}
 \author{Tom Kennedy}
 \homepage{http://www.math.arizona.edu/~tgk}
 \email{tgk@math.arizona.edu}
\affiliation{Departments of Mathematics and Physics 
\\ University of Arizona, Tucson, AZ, 85721}

\date{\today}

\begin{abstract}
The conjecture that the scaling limit of the two-dimensional 
self-avoiding walk (SAW)
in a half plane is given by the stochastic Loewner evolution (SLE) 
with $\kappa=8/3$ leads to explicit predictions about the SAW. 
A remarkable feature of these predictions
is that they yield not just critical exponents, but probability
distributions for certain random variables associated with the 
self-avoiding walk. We test two of these predictions with 
Monte Carlo simulations and find excellent agreement, thus providing 
numerical support to the conjecture that the scaling limit of the SAW
is SLE$_{8/3}$.
\end{abstract}

\maketitle

A variety of two dimensional models in statistical physics
are expected to have conformally invariant scaling limits.
This conformal invariance has made it possible
to determine critical exponents for these two dimensional models.
Recently, Schramm introduced a two dimensional conformally invariant 
random process which he called stochastic 
Loewner evolution \cite{schramm}. 
This process depends on a parameter $\kappa$, and so is 
denoted SLE$_{\kappa}$. The qualitative behavior of the process 
changes with the parameter and it appears that for different 
values of $\kappa$, the process is related to the scaling
limit of various two dimensional models. 

Schramm's SLE process is defined in a half plane as follows. 
Let $B_t$ be a standard
one-dimensional Brownian motion. We define a complex-valued function 
$g_t(z)$ for $z$ in the upper half of the complex plane by 
the following differential equation,
\begin{equation} 
\partial_t g_t(z) = {2 \over g_t(z) - \sqrt{\kappa} B_t},
\end{equation} 
along with the initial condition $g_0(z)=z$.
The SLE trace is the curve defined by 
\begin{equation}
\gamma(t)=\lim_{z \rightarrow 0} g_t^{-1}(z+\sqrt{\kappa} B_t).
\end{equation}
(The limit is only over $z$ in the upper half plane. 
$g_t^{-1}$ can be obtained by solving the differential equation backwards 
in time.)
The limit exists and gives a continuous curve. This has been proved for 
$\kappa \ne 8$, and is believed to be true for all $\kappa$ \cite{schramm,rs}.
For $\kappa=8/3$ it is conjectured that this SLE trace gives the scaling 
limit of the SAW restricted to the half plane \cite{lswconj}. 

Schramm showed that if the loop-erased random walk has a conformally invariant
scaling limit, then that limit must be SLE$_2$ \cite{schramm}. 
He also conjectured that the scaling limit of percolation 
should be related to SLE$_6$, and the scaling limit of uniform spanning
trees is described by SLE$_2$ and SLE$_8$.
Smirnov has proved the conformal invariance conjecture
for critical percolation on the triangular lattice and that
SLE$_6$ describes the limit \cite{smirnov}.
Lawler, Schramm and Werner used SLE$_6$ to rigorously determine the 
``intersection exponents'' for Brownian motion 
and proved a conjecture of Mandelbrot that the outer boundary of a Brownian
path has Hausdorff dimension 4/3 \cite{lswbm, lswa, lswb, lswc}.
Rohde and Schramm conjectured that the random cluster representation
of the Potts model for $0<q<4$ is related to the SLE process as well
\cite{rs}.

\begin{figure}
\includegraphics{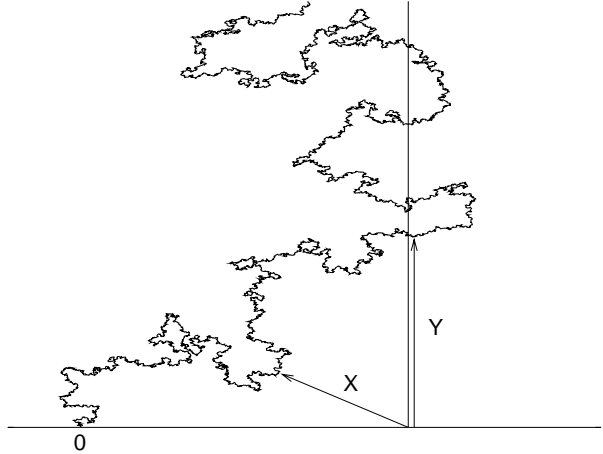}
\caption{\label{fig:rv} The random variables $X$ and $Y$ that we study
are illustrated for a SAW. (Even though it appears the 
path intersects itself, when viewed at a smaller resolution  one sees 
that it does not.)
}
\end{figure}

By using stochastic calculus it is possible to compute many quantities
related to SLE. A remarkable feature of these calculations is that 
they can yield not just critical exponents, 
but entire probability distributions for certain random variables. 
In this paper we will consider two examples of such random variables,
illustrated in figure \ref{fig:rv}.
For SLE$_{8/3}$ one can explicitly compute the probability distributions
of $X$ and $Y$. We will compare these explicit distributions
with Monte Carlo simulations of the random variables for the SAW.

To define the first random variable we 
fix a point $(c,0)$ on the horizontal axis. Given an SLE trace or a SAW
with lattice spacing $\delta$,
we consider the distance from the curve to the point $(c,0)$.
We define $X$ to be the ratio of this distance to $c$. 
\begin{equation} 
X={1 \over c}  \min_{t \ge 0}  ||\gamma(t)-(c,0)|| 
\end{equation} 
So $X$ takes values in $(0,1]$. SLE is invariant under dilations, 
so the distribution of $X$ is independent of $c$. In the scaling limit
($\delta \rightarrow 0$) the distribution of $X$ for the SAW should also
be independent of $c$. 

For the second random variable 
we consider the intersections of the curve with the vertical line $x=c$. 
We define $Y$ to be the ratio of the $y$-coordinate of the
lowest intersection to $c$. 
\begin{equation} 
Y={1 \over c} \min \{y: (y,c)=\gamma(t) \,\, for \,\, some \,\, t \ge 0\} 
\end{equation} 
So $Y$ takes values in $(0,\infty)$. The distribution of $Y$ should also 
be independent of $c$. 

The distributions of $X$ and $Y$ will follow from a 
remarkable theorem of Lawler, Schramm and Werner \cite{lswconj} 
for SLE$_{8/3}$, 
which appears in a survey article by Lawler 
\cite{lawler} about their joint results.
Let $H$ be the upper half plane, and $\bar H$ its closure. 
Let $A$ be a compact subset of $\bar H$ which does not contain $0$ and such 
that $H \setminus A$ is simply connected.
Let $\Phi_A$ be the conformal map from $H \setminus A$ onto $H$ which 
fixes $0$ and $\infty$ and has $\Phi_A^\prime(\infty)=1$.

\smallskip

\noindent {\it Theorem (Lawler, Schramm, Werner) 
For $\kappa=8/3$, SLE in a half plane satisfies
\begin{equation}
P(\gamma[0,\infty) \cap A = \emptyset) = \Phi^\prime_A(0)^{5/8} 
\end{equation}
}

\smallskip

The Riemann mapping theorem says that the conformal map $\Phi_A$ 
exists, but an explicit formula for it is available in only a few 
special cases. For the two random variables 
defined above, the computation of their distributions is just an application 
of the above theorem for two $A$'s for which there is an explicit 
conformal map. 

We start with the random variable $X$, and take $c=1$. 
For $a<1$, let $A_a$ be the half of the disc centered at $(1,0)$ 
with radius $a$ that is in the upper half plane. 
The distance $X$ from $\gamma[0,\infty)$ to $(1,0)$ is less than or equal 
to $a$ if and only if $\gamma[0,\infty)$ hits $A_a$. So 
\begin{equation} 
P(X \le a) = P(\gamma[0,\infty) \cap A_a \ne \emptyset) 
\end{equation}

\begin{figure}
\includegraphics{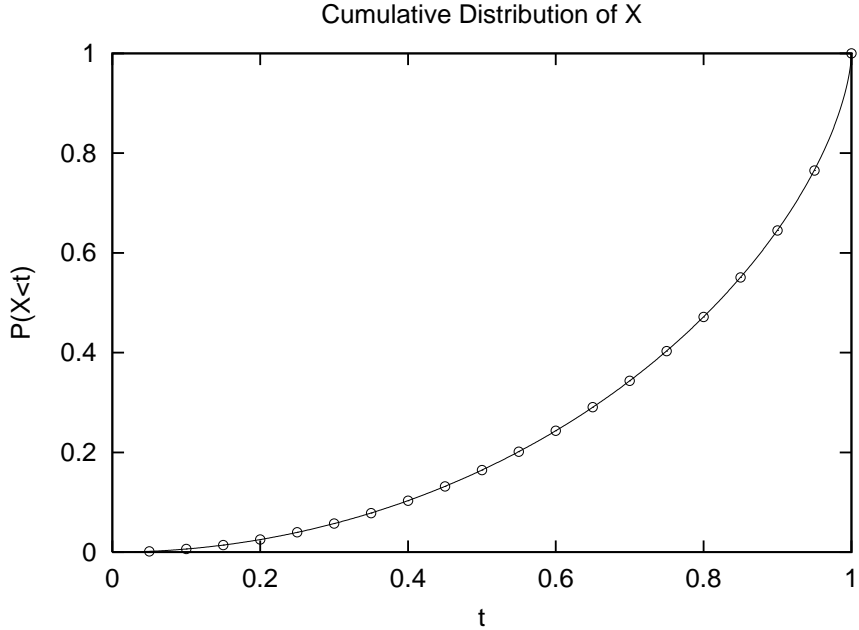}
\caption{\label{fig:dist_exact} The curve is the exact distribution of $X$ 
for SLE$_{8/3}$. The points are the results of a 
Monte Carlo simulation for the SAW on the square lattice.
}
\end{figure}

The conformal map that sends $H \setminus A_a$ onto $H$ is 
\begin{equation} 
\Phi_{A_a}(z) = z-1 + {a^2 \over z-1} + 1 + a^2,
\end{equation}
It is normalized so that it fixes $0$ and $\infty$ and 
has $\Phi_{A_a}^\prime(\infty)=1$. 
Since $\Phi^\prime_{A_a}(0) = 1-a^2$,
the theorem says 
\begin{equation} 
P(X \le t) = 1-(1-t^2)^{5/8}.
\end{equation}
In figure \ref{fig:dist_exact} the solid curve is the above function.
The results of our Monte Carlo simulation for the SAW on the square lattice
are shown with circles. One cannot see any difference at the scale of 
this plot. 

Now we find the distribution of $Y$. Again, we take $c=1$. 
Let $L_a$ to be the vertical line
segment from $(1,0)$ to $(1,a)$. 
The random variable $Y$ is less than or equal to $a$ 
if and only if the curve hits $L_a$. So 
\begin{equation} 
P(Y \le a) = P(\gamma[0,\infty) \cap L_a \ne \emptyset) 
\end{equation}

\begin{figure}
\includegraphics{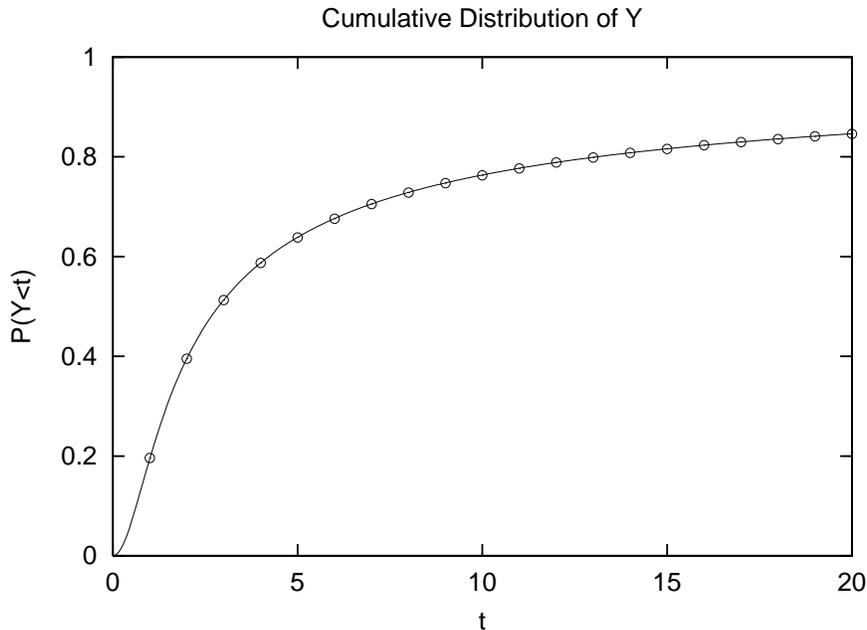}
\caption{\label{fig:vert_exact} The curve is the exact distribution of $Y$ 
for SLE$_{8/3}$. The points are the results of a 
Monte Carlo simulation for the SAW on the square lattice.
}
\end{figure}

The conformal map that takes $H \setminus L_a$ onto the upper half plane is 
\begin{equation}
\Phi_{L_a}(z)= i \sqrt{-(z-1)^2 -a^2} + \sqrt{1+a^2}
\end{equation}
(The square root is defined with the 
usual branch cut along the negative $x$-axis.) 
The map fixes $0$ and $\infty$ and has $\Phi_{L_a}^\prime(\infty)=1$. 
We have
$\Phi_{L_a}^\prime(0) = 1 / \sqrt{1+a^2}$,
so the theorem says 
\begin{equation}
P(Y \le t) = 1-(1+t^2)^{-5/16}
\end{equation}
In figure \ref{fig:vert_exact} the solid curve is the above function and 
the results of our Monte Carlo simulation for the SAW on the square lattice
are again shown with circles. 

We simulate the self-avoiding walk using the pivot algorithm 
(Expository accounts of this algorithm 
may be found in \cite{madrasslade,sokala}.)
An iteration of this Markov chain algorithm
starts by picking a random site on the walk. 
Then one picks a random lattice symmetry $g$.
The section of the walk from the starting point to the randomly chosen 
site is not changed. The rest of the 
walk is ``pivoted'' by applying $g$ to it with respect to the randomly 
chosen site. If the resulting walk is self-avoiding, it is accepted;
otherwise it is rejected. 
The pivot algorithm can also be used to simulate the SAW restricted 
to the half plane. When a pivot is proposed, in addition to 
checking if the pivoted walk is self-avoiding, one must also check 
if it stays in the half-plane. If both of these conditions are satisfied,
then the pivot is accepted.  

If one uses a hash table and checks for self-intersections in the 
pivoted walk starting from the pivot point and 
working outwards, then the time required
for the pivot algorithm to produce an accepted pivot is believed to be
$O(N)$ on average \cite{madrassokal}. 
It has been shown recently that by taking advantage of the nearest neighbor
nature of the walk when checking for self-intersections and using a 
data structure to store the walk that postpones carrying out the pivots, 
the pivot algorithm in two dimensions can be implemented so that the 
time required to produce an accepted pivot is $O(N^q)$ with $q$ less 
than one \cite{tk_pivot}. The exact value of $q$ is not known, but it 
appears to be less than $0.57$ for two-dimensional walks. 

\begin{figure}
\includegraphics{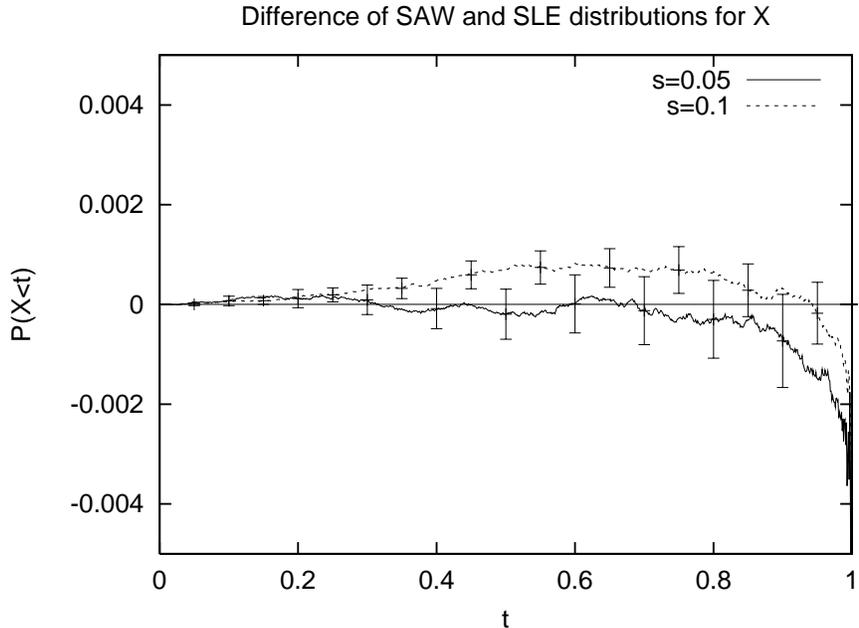}
\caption{\label{fig:dist} Each curve is the difference between the 
Monte Carlo computation of the cumulative distribution of $X$ for the 
SAW and the exact cumulative distribution of $X$ for SLE$_{8/3}$. 
Error bars of two standard deviations are shown for selected values of $t$. 
}
\end{figure}

To study the scaling limit of the SAW walk there are two limits that
must be taken. We must let the number of steps, $N$, go to infinity
and we must take the lattice spacing to zero. A walk with $N$ steps 
lives on scale $N^{3/4}$, so to study the random variables $X$ and $Y$ 
we take $c=s N^{3/4}$. 
We must take $s$ small to make the effect of the finite length of our 
walks negligible, and we must take $N$ large to make the lattice spacing small
compared to the scale of the random variables.
Our simulations are done on the square lattice with $N=1,000,000$. 
We study the distribution of $X$ for $s=0.05$ and $s=0.1$. 
For $Y$ we use $s=0.005$ and $s=0.01$. 
We compute $P(X \le t)$ and $P(Y \le t)$ for $1,000$ values of $t$. 
For $X$, the values of $t$ range from $0$ to $1$. 
For $Y$ the values range from $0$ to $20$. 

There are 20 billion iterations of the Markov chain in the simulation.
Approximately $5\%$ of these proposed pivots are accepted. 
For some random variables, e.g., the end to end distance of the walk, 
an accepted pivot produces a radical change in the value of the random
variable, and it is expected that each accepted pivot produces an 
essentially independent sample of the random variable. For our
random variables a pivot can change their values 
only if it occurs in a relatively short segment of the 
walk near the origin. So the autocorrelation time will be significantly 
longer than for random variables like the end to end distance. 

\begin{figure}
\includegraphics{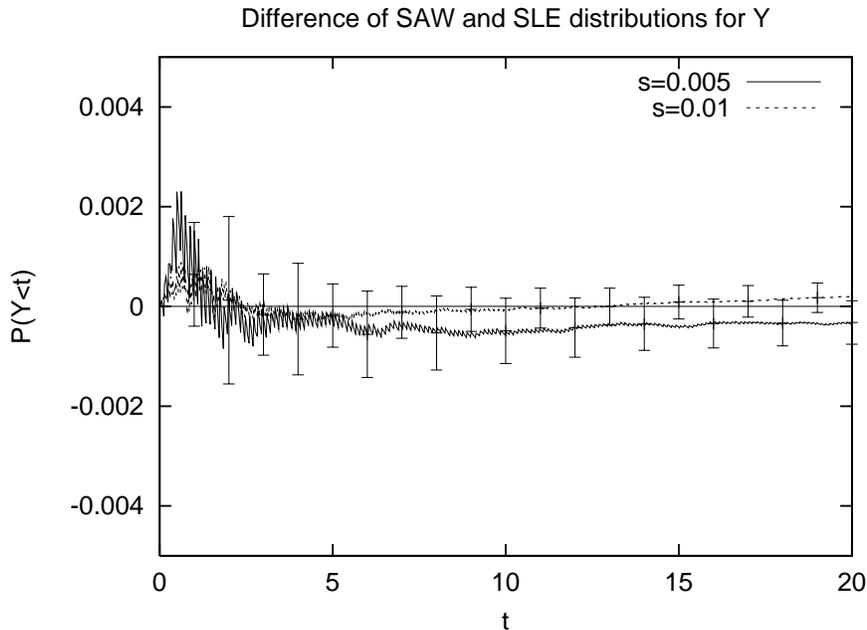}
\caption{\label{fig:vert} Each curve is the difference between the 
Monte Carlo computation of the cumulative distribution of $Y$ for the 
SAW and the exact cumulative distribution of $Y$ for SLE$_{8/3}$. 
Error bars of two standard deviations are shown for selected values of $t$. 
}
\end{figure}

Since it is impossible to see the difference between the SAW simulations
and the SLE curves in figures \ref{fig:dist_exact} and \ref{fig:vert_exact}, 
we plot the distributions from the SAW simulations minus
the SLE$_{8/3}$ distributions in figures 
\ref{fig:dist} and  \ref{fig:vert}. The most important feature of 
these figures is the scale on the vertical axis. It is $1/100$ of 
the scale of figures \ref{fig:dist_exact} and \ref{fig:vert_exact}. 
So the agreement between the SLE distributions and the simulated 
SAW distributions is excellent. 
Since the simulation uses finite length walks with a nonzero lattice
spacing, if the simulation was run long enough we would see that 
the quantity being plotted is not exactly zero. 
Error bars, based on two standard deviations, are shown in the figures
for selected values of $t$. (Different values of 
$t$ are used for the two curves so the error bars do not overlap.)
The error bars are the same order of magnitude as the quantity 
being plotted, so our simulation is not quite accurate enough to 
see the finite length and lattice effects for the values of $N$ and $s$
that we use. (In the $s=0.1$ curve for $X$, the effects might be just 
beginning to emerge.)

Our simulations of the SAW walk in a half plane have shown that 
the distributions of two particular random variables related to the 
walk agree extremely well with the exact distributions of SLE$_{8/3}$
for these random variables. This supports the conjecture that the 
scaling limit of the SAW is SLE$_{8/3}$. 
Schramm has recently 
given a formula for the probability that the SLE curve passes to 
the right of a fixed point in the half plane for all values of 
$\kappa$ \cite{schramm_perc}.
We expect that the distributions
of many more random variables associated with SLE$_{8/3}$ will be 
found in the near future. 

This work was supported by the National Science Foundation (DMS-9970608).

\bibliography{saw_sle}

\end{document}